\documentstyle{amlts}
\input classical12.sty
\begin{document}
\annalsline{155}{2002}
\received{March 11, 2000}
\startingpage{131}
\def\bye{\end{document}}
 \font\tenrm=cmr10

\input amssym.def
\input amssym.tex
\input boxedeps.tex 
\SetepsfEPSFSpecial 
\HideDisplacementBoxes
\def\figin#1#2{ 
$$
 {\BoxedEPSF{#1.eps scaled
#2}%
}%
$$
}

\def\ritem#1{\item[{\rm #1}]}
\def\End{{\rm End}}
\def\be{\begin{equation}}
\def\ee{\end{equation}}
\def\barr#1{\overline{#1}}
\def\Sym{{\rm Sym}}
\def\ra{\rightarrow}
\def\GL{{\rm GL}}
\def\zz{{\Bbb Z}}
\def\Aut{{\rm Aut}}
\def\rr{{\Bbb R}}
\def\cc{{\Bbb C}}
\catcode`\@=11
\font\twelvemsb=msbm10 scaled 1100
\font\tenmsb=msbm10
\font\ninemsb=msbm10 scaled 800
\newfam\msbfam
\textfont\msbfam=\twelvemsb  \scriptfont\msbfam=\ninemsb
  \scriptscriptfont\msbfam=\ninemsb
\def\msb@{\hexnumber@\msbfam}
\def\Bbb{\relax\ifmmode\let\next\Bbb@\else
 \def\next{\errmessage{Use \string\Bbb\space only in math
mode}}\fi\next}
\def\Bbb@#1{{\Bbb@@{#1}}}
\def\Bbb@@#1{\fam\msbfam#1}
\catcode`\@=12

 \catcode`\@=11
\font\twelveeuf=eufm10 scaled 1100
\font\teneuf=eufm10
\font\nineeuf=eufm7 scaled 1100
\newfam\euffam
\textfont\euffam=\twelveeuf  \scriptfont\euffam=\teneuf
  \scriptscriptfont\euffam=\nineeuf
\def\euf@{\hexnumber@\euffam}
\def\frak{\relax\ifmmode\let\next\frak@\else
 \def\next{\errmessage{Use \string\frak\space only in math
mode}}\fi\next}
\def\frak@#1{{\frak@@{#1}}}
\def\frak@@#1{\fam\euffam#1}
\catcode`\@=12

 \def\co{\colon\,}
\newcommand{\Inv}{L}   
\newcommand{\ab}{\mbox{\scriptsize ab}}
\newcommand{\old}{\mbox{\scriptsize old}}
\newcommand{\Ref}{{\rm Ref}}
\newcommand{\LF}{{\rm LF}}
\newcommand{\HP}{{\rm HP}}
\newcommand{\TP}{{\rm TP}}
\newcommand{\Pro}{{\rm Pro}}
\newcommand{\Brd}{{\rm GB}}
\newcommand{\GB}{{\rm GB}}

\title{Braid groups are linear}

 \acknowledgements{Financially supported by the Swiss National Science Foundation. }
 \author{Daan Krammer}
 \institutions{University of Basel, Basel, Switzerland\\
{\eightpoint {\it E-mail address\/}: Daan.Krammer@unibas.ch}}

\bigbreak\centerline{\bf Abstract}
\bigbreak

In a previous work \cite{kra}, the author considered a representation of the braid group 
$\rho\co B_n\ra\GL_m(\zz[q^{\pm1},t^{\pm1}])$ ($m=n(n-1)/2$), and proved it to be faithful for $n=4$. Bigelow
\cite{big2} then proved the same representation to be faithful for all $n$ by a beautiful topological argument. The
present paper gives a different proof of the faithfulness for all $n$. We establish a relation between the Charney
length in the braid group and exponents of $t$. A certain $B_n$-invariant subset of the module is constructed whose
properties resemble those of convex cones. We relate line segments in this set with the Thurston normal form of a
braid.

\medbreak \centerline{\bf Contents}

1.  Introduction

2. Combinatorial preliminaries

3. The representation

4.  Faithfulness

5.  Half-permutations

6.  Two more properties of the representation

\hglue14pt References

\section{Introduction} \label{seca}
\demo{Statement and history of the problem} 
A group is said to be linear if it is isomorphic to a subgroup of $\GL(n,K)$ for some natural number $n$ and some field $K$. An interesting question asks whether the braid group is linear.

One of the most famous representations of the braid group
 is the Burau representation $B_n\ra\GL_{n-1}(\zz[q^{\pm1}])$. It is easily shown to be faithful for $n\leq3$. Moody
[15] proved the Burau representation to be unfaithful for $n\geq9$. This bound was improved to $n\geq6$ by
Long and Paton [13] and to $n\geq5$ by Bigelow \cite{big}. It is still unknown whether the Burau
representation of $B_4$ is faithful.

One of the braid group representations, previously studied 
by Lawrence [12], was proved to be faithful by the author [11] in the case of $B_4$. Shortly
thereafter, Bigelow \cite{big2} found a proof that the same representation is faithful for all $n$ by a
beautiful topological argument. The present article deals with again the same representation.

More on the history of the linearity problem for braid groups can be found in Birman's review \cite{bir2}.
\enddemo

\demo{The representation} 
The representation of our interest will be denoted $\rho\co B_n\ra\GL(V)$, where
$V$ is an $m$-dimensional free module over some ring $R$, with $m=n(n-1)/2$.
 It depends on two invertible elements $q,t\in R$. There are many definitions of this representation. This paper
follows an elementary route by exhibiting the entries of the involved matrices, and completely avoids topological
arguments. Other definitions include a second homology group \cite{kra}, \cite{law} and a pictorial approach
\cite{big2}, \cite{kra}. Zinno \cite{zin} recently showed the representation to be a summand of the Birman-Wenzl
algebra \cite{bw}.
\enddemo

\demo{Combinatorial preliminaries} 
Our linearity proof for the braid group involves a solution to the word
problem in the braid group. Among the solutions to the word problem we
mention Artin's one \cite{art} ($B_n$ is isomorphic to a subgroup of
$\Aut(F_n)$) and a solution based on Thurston's boundary of Teichm\"uller space \cite{pen}. 
Neither solution is relevant to this paper. Important for us is a third, again totally different solution due to Garside
(\cite{gar}, see also \cite{del}, \cite{C+}).

For $1\leq i<j\leq n$, let $s(i,j)=s_{ij}$ denote the permutation (called a reflection) 
in the symmetric group $S_n$ interchanging $i$ with $j$ and preserving the rest. The set of reflections in $S_n$ will
be denoted by $\,\Ref$. Let $\ell\co S_n\ra\zz_{\geq0}$ denote the length function with respect to
$\{s_{12},s_{23},\ldots,s_{n-1,n}\}\subset S_n$.

The braid group $B_n$ admits a presentation by generators $\{rx\mid x\in S_n\}$
 and relations $r(xy)=(rx)(ry)$ whenever $\ell(xy)=\ell(x)+\ell(y)$. The positive braid monoid $B_n^+$ is by definition
the submonoid of $B_n$ generated by $\Omega:=r(S_n)$. For $x\in B_n^+$ there exists a unique longest $y\in S_n$ with
$x\in(ry)B_n^+$, notation:  $ry=\LF(x)$. We will make use of the following proposition, which is implied by Garside's
results.
\enddemo

\proclaimtitle{See \ref{fa28}}
\specialnumber{A}
\proclaim{Proposition}  Let $B_n$ act on a set $U${\rm .} Suppose we are given nonempty
disjoint subsets $C_x\subset U$ ($x\in\Omega$) with $xC_y\subset C_{\LF(xy)}$ for all $x,y\in\Omega${\rm .} Then the
$B_n$\/{\rm -}\/action on $U$ is faithful.\endproclaim

Later on, we will apply Proposition~A by putting $U=V$. The central question is to find $C_x$ satisfying the
assumption of Proposition~A. The $C_x$ we use will be convex in some sense.
\pagebreak

For any $x\in B_n^+$ there is a unique $(x_1,\ldots,x_k)\in\Omega^k$ such that $x_1\cdots x_k=x$ and
$\LF(x_ix_{i+1})=x_i$ for all $i$, and $x_k=1$. It is called the greedy form of $x$ and is due to Garside \cite{gar}.
Thurston \cite{C+} showed that any braid $x\in B_n$ can uniquely be written $x=y^{-1}z$ with $y,z\in B_n^+$ such that
there is no $w\in B_n^+ -\{1\}$ with $\{y,z\}\subset wB_n^+$. Writing $(y_1,\ldots,y_k)$ for the greedy form for $y$
and $(z_1\,\ldots,z_\ell)$ for the greedy form for $z$ then gives $x=y_k^{-1}\cdots y_1^{-1}z_1\cdots z_\ell$; this
is called the Thurston normal form. Closely related is the length function $\ell_\Omega\co B_n\ra\zz_{\geq0}$ with
respect to $\Omega$. Charney \cite{cha} showed that the growth function $\sum_{x\in
B_n}z^{\ell_\Omega(x)}\in\zz[[z]]$ is rational. We call $\ell_\Omega$ the Charney length function.

\demo{Faithfulness} 
We will throughout make use of a certain basis $\{x_s\mid s\in\Ref\}$ of $V$, and will identify an element of $\End(V)$ with its matrix with respect to this basis. Thus, $\End(V)$ is identified with $M_m(R)$, the size $m$ matrix algebra over $R$.

We will observe that $\rho B_n^+\subset M_m(\zz[q,q^{-1},t])$; i.e., for positive braids $x$, the entries of $\rho x$ do
not involve negative powers of $t$. 

Henceforth, we assume $R=\rr[t^{\pm1}]$, $q\in\rr$ and $0<q<1$. Then for all positive braids $x\in B_n^+$, the entries of $\rho x$ are in $\rr_{\geq0}+t\,\rr[t]$. This observation is the most important step of the faithfulness proof of $\rho$. A faithfulness proof of the braid group seems to be impossible without some kind of inequalities involved (think of convex cones), and the foregoing observation fulfills this need.

Let $M_m(\{0,1\})$ denote the set of size $m$ square matrices with entries in $\{0,1\}\subset\zz$. Multiplication in $M_m(\{0,1\})$ is defined as follows. Given two elements, one first multiplies them in $M_m(\zz)$, then replaces all positive entries by one, leaving zero entries untouched. This multiplication makes $M_m(\{0,1\})$ into a monoid. We have a monoid homomorphism $B_n^+\ra M_m(\{0,1\})$, the image of $x\in B_n^+$ being obtained from $\rho x$ by setting $t=0$ and then replacing positive entries by one. Now $M_m(\{0,1\})$ is finite; the combinatorics of the homomorphism $B_n^+\ra M_m(\{0,1\})$ are crucial in the correct definition of $C_x$, which is briefly as follows.

Define 
\begin{eqnarray*}
 \HP&=&\Big\{A\subset\Ref\ {\Big\vert}\ s_{ij},s_{jk}\in A\Rightarrow s_{ik}\in A\mbox{ whenever
}1\leq i<j<k\leq n\Big\}, \\
 \Inv(x)&=&\Big\{s_{ij}\ {\Big\vert}\ 1\leq i<j\leq n,\ x^{-1}i>x^{-1}j\Big\},\ \ \ (x\in S_n).
\end{eqnarray*}
We will see that for any $A\in\HP$ there is a (unique) greatest $B\in\Inv(S_n)$ with $B\subset A$.
Notation: 
$B=\Pro(A)$. For $x\in\Omega$, one defines $C_x\subset V$ to be the set of those vectors $\sum_{s\in\Ref}a_sx_s$
with $a_s\in\rr_{\geq0}+t\,\rr[t]$ and such that on putting $A:=\{s\in\Ref\mid a_s\in t\,\rr[t]\}$ one has $A\in \HP$
and $x=r\,\Inv^{-1}\,\Pro(A)$. 

Clearly, it is a purely combinatorial issue whether $xC_y\subset C_{\LF(xy)}$ for all $x,y\in\Omega$ (the condition of Proposition~A). It turns out to be correct, whence by Proposition~A:

\proclaimtitle{See \ref{tf7}}
\specialnumber{B}
\proclaim{Theorem}
 The representation $\rho\co B_n\ra\GL(V)$ is faithful{\rm ,} even if $q$ is a real number with
$0<q<1${\rm .}\endproclaim

 Theorems~C and D below state two closely related properties of the
representation. They are new and will be proved in Section~\ref{secf}.

\proclaimtitle{See \ref{tf26}}
\specialnumber{C}
\proclaim{Theorem} Let $x\in B_n${\rm ,} and consider the Laurent expansion of $\rho x$ with
respect to $t${\rm :}
$$ \rho x=\sum_{i=k}^\ell A_i(q)\,t^i,\ \ \ A_i\in M_m(\zz[q^{\pm1}]),\ \ \ A_k\neq0,\ \ \ A_\ell\neq0. $$
\smallbreak
{\rm (a)} Then $\ell_\Omega(x)=\max(\ell-k,\ell,-k)${\rm .}

\smallbreak {\rm (b)} If{\rm ,} in addition{\rm ,} $x\in B_n^+ -\Delta B_n^+${\rm ,} then $k=0$ and
$\ell=\ell_\Omega(x)${\rm .} 
\endproclaim
 
We define an ordering on $R=\rr[t^{\pm1}]$ by giving a nonzero element of it~the same sign as its trailing
coefficient (the coefficient for the least occurring exponent of $t$). We write $C_1=C$ and
$U=\mathbold{\cup}_{x\in B_n}xC$. The following result shows that $U$ has properties resembling those of
convex cones in real vector spaces, and moreover connects the Thurston normal form with line segments in
$U$.

\vglue-16pt
\phantom{try}

\proclaimtitle{See \ref{tf28}}
\specialnumber{D}
\proclaim{Theorem}  \vglue-20pt
\phantom{trick}
  \begin{itemize}
\ritem{(a)} The $xC$ (with $x\in B_n$) are disjoint{\rm .}
\ritem{(b)} Let $(\tilde{y}_1,\ldots,\tilde{y}_k)$ be a Thurston normal form\/{\rm ;} i.e.{\rm ,} there are
 greedy\break
$(u_1,\ldots,u_s)${\rm ,}
$(v_1,\ldots,v_t)$ with $(u_s^{-1},\ldots,u_1^{-1},v_1,\ldots,v_t)=(\tilde{y}_1,\ldots,\tilde{y}_k)${\rm ,}
 and $u_s,v_t\neq
1${\rm ,} and there is no $w\in B_n^+ -\{1\}$ such that $\{u_1,v_1\}\subset wB_n^+$. Let $\tilde{x}_0,\ldots,
\tilde{x}_k\in
B_n$ be such that $\tilde{x}_i=\tilde{x}_{i-1}\tilde{y}_i$ $(1\leq i\leq k).$ Then 
$$ \frac{t^i\tilde{x}_0C+\tilde{x}_kC}{t^i+1}\subset \left\{ \begin{array}{llr}
{\tilde{x}_0}C,\ \ & i\leq -s; \\
{\tilde{x}_{i+s}}C, & -s\leq i\leq t; \\
{\tilde{x}_k}C, & t\leq i. \\ \end{array} \right. $$
\ritem{(c)} The set $U$ is closed under addition and scalar multiplication by positive elements of $R${\rm .}
\end{itemize}
\phantom{trick}
\vglue-16pt

\endproclaim

{\it Comparison of three methods}.
In a previous paper \cite{kra}, the representation $\rho$ is proved to be faithful for $n=4$ by a somewhat different 
method. I do not know whether this method works for $n>4$. The differences and similarities between this method and
the method of the present paper are as follows. Briefly, the roles (not the meanings) of $q$ and $t$ are interchanged. 

One of our results, Theorem~C, relates the exponents of $t$ with the Charney length function. In \cite{kra} one finds a (for $n>4$ conjectural) relation between the exponents of $q$ and the length function with respect to some other generating subset $Q\subset B_n$ with cardinality 
$$ |Q|=\frac{1}{n+1}{{2n \choose n}}. $$
A basic reference to $Q$, which is also known as the set of band generators, is~\cite{bir1}. The present paper
assumes $q$ to be a real number with $0<q<1$; in \cite{kra}, $t$ is a real number with $0<t<1$.

The present paper studies the set
$$ \bigoplus_{s\in\Ref}\Bigl( \rr_{\geq0}+t\,\rr[t] \Bigr) \,x_s\subset V, $$
which is essentially a simplicial cone. If $t=1$, then the $B_n$-module $V$ can be shown to be the symmetric square of the Burau module, so that `the cone of positive semi-definite elements' makes sense. In \cite{kra}, a generalization of the cone of positive semi-definite elements is studied. This convex cone is not simplicial at all; rather, it is given by finitely many nonlinear algebraic inequalities.

A third method of proof was found by Bigelow \cite{big2}. 
His beautiful and strikingly short proof involves neither a solution to the word problem, nor a basis of the module. In
Bigelow's proof, both $q$ and $t$ are variables. The total ordering on $\langle q,t\rangle $ he uses makes $q$ ``more important''
than $t$, so that his method is closer to having $t$ constant than to having $q$ constant.

It seems to be interesting to combine the three approaches into one theory, which presumably involves both generating sets $Q$ and $\Omega$.

\demo{Overview} 
The paper is built as follows. There are two sets of combinatorial
results. The first set is mainly due to Garside, Thurston and Charney
and is collected in Section~\ref{secb}. The second set might be new and
is treated in Section~\ref{sece}. In Section~\ref{secc}, we define the
representation and establish a few identities. An overview of the
faithfulness proof (but more detailed than in the introduction) can be
found in Section~\ref{secd}. Section~\ref{secf} is devoted to proving Theorems~C and D.

\demo{Acknowledgements}
 The author gratefully acknowledges the support by the Swiss National Science Foundation.
Many thanks to Hanspeter Kraft for his warm support. \pagebreak

\section{Combinatorial preliminaries} \label{secb}

This section collects some combinatorial properties of braid groups mainly due to Garside, Thurston and Charney. For proofs, we refer to \cite{gar}, \cite{del}, \cite{C+}, \cite{cha}, \cite{mic}; remaining statements are left to the reader to prove.

 The braid group $B_n$ is defined to be the fundamental group of $\{X\subset\cc\co |X|=n\}$, the set of $n$-element
subsets of $\cc$, with its obvious topology. Artin proved that the braid group $B_n$ admits a finite presentation
(called the {\it Artin presentation}) with generators $\sigma_1$,\ldots,$\sigma_{n-1}$ and relations
\begin{eqnarray}
 \sigma_i\sigma_{i+1}\sigma_i&=&\sigma_{i+1}\sigma_i\sigma_{i+1},\ \ \ (1\leq i\leq n-1),
\label{tf1}
\\[6pt]
\sigma_i\sigma_j&=&\sigma_j\sigma_i,\ \ \ (|i-j|>1). \label{tf2} \end{eqnarray}
(We will view $\sigma_i$ as an element of the braid group.)  

Let $S_n$ denote the symmetric group on $I_n=\{1,2,\ldots,n\}$ 
(action from the left). For $1\leq i<j\leq n$, let $s_{ij}=s(i,j)\in S_n$ denote the permutation (called a {\it reflection}) interchanging $i$ with $j$ and fixing the other elements of $I_n$. Put $s_i=s_{i,i+1}$ and
$S=\{s_1,\ldots,s_{n-1}\}$. (The pair $(S_n,S)$ is known as a Coxeter system of type $A_{n-1}$.) By $\Ref$ we will
denote the set of reflections in $S_n$.

Let $\ell\co S_n\ra\zz_{\geq0}$ denote the length function with respect to $S$; i.e., $\ell(x)$ is the smallest natural
number $k$ such that there exist $x_1,\ldots,x_k\in S$ with $x=x_1\cdots x_k$. The symmetric group $S_n$ contains
a unique longest element $w_0$, given by $w_0(i)=n+1-i$.

The braid group $B_n$ admits a presentation with generators $\{rx\mid x\in S_n\}$ 
and relations $r(xy)=(rx)(ry)$ whenever $\ell(xy)=\ell(x)+\ell(y)$. We will view $rx$ as an element of $B_n$, and we
denote the image of $r\co S_n\ra B_n$ by $\Omega$. There exists a well-known homomorphism $B_n\ra S_n$ defined
by $rx\mapsto x$ $(x\in S_n)$. One can identify $r(s_i)$ with $\sigma_i$ in the Artin presentation of the braid group.
The element $\Delta:=r(w_0)$ is known as the {\it half-twist}.

The submonoid of $B_n$ generated by $\Omega$ will be denoted $B_n^+$ (this includes~1). 
Elements of the braid group $B_n$ are called {\it braids\,} and elements of $B_n^+$ are called {\it positive braids}.
Recall the length function $\ell\co S_n\ra\zz_{\geq0}$. By the same symbol, we will denote the {\it length
function\,} $\ell\co B_n^+\ra\zz_{\geq0}$, which is the (unique) monoid homomorphism with $\ell(rx)=\ell(x)$ for all
$x\in S_n$. Let $\Omega_k$ denote the set of elements of $\Omega$ of length $k$.

A {\it smallest\,} (respectively, {\it greatest}) element of a (partially)
 ordered set is an element which is smaller (respectively, greater) than any other element. A smallest or greatest
element does not necessarily exist, but if it exists, then it is unique.

Define an ordering on $B_n^+$ by $x\leq y\Leftrightarrow y\in xB_n^+$. 
Restriction of this ordering yields an ordering on $\Omega$, and thereby on $S_n$. The ordering on $S_n$ can
equivalently be given by $x\leq xy$ if and only if $\ell(xy)=\ell(x)+\ell(y)$; it is known as the weak Bruhat ordering.
The ordered set $S_n$ has a smallest element $1$ and a greatest element $w_0$. The smallest element of $\Omega$
is also denoted 1, and its greatest element is $\Delta$.

It can be shown that for any $x\in B_n^+$, the set $\{y\in \Omega\mid y\leq x\}$
 has a greatest element. It will be denoted by $\LF(x)$ (Left most Factor). A sequence $(x_1,\ldots,x_k)\in \Omega^k$
is said to be  (\/{\it left\/{\rm )} greedy\,} if $\LF(x_ix_{i+1})=x_i$ for all $i=1,\ldots,k-1$. For any $x\in B_n^+$,
there is a unique greedy sequence $(x_1,\ldots,x_k)$ with $x_1\cdots x_k=x$ and $x_k\neq 1$. It is called the  
(\/{\it left\/{\rm )}  greedy form\,} for $x$.

An important identity reads
\be \LF(xy)=\LF(x\,\LF(y)) \label{tf4} \ee 
for all $x,y\in B_n^+$. It implies that the map $B_n^+\times \Omega\ra \Omega$ defined by $(x,y)\mapsto\LF(xy)$ is an action of the monoid $B_n^+$ on $\Omega$.

The following proposition singles out an aspect of the
 word problem which will be used in the present paper. A similar result can be found in \cite{kra}. Its proof is a
simple application of Garside's results described above. The result gives a sufficient condition on a $B_n$-action on
any set to be faithful. Later on, the set will be chosen to be a module.

\proclaim{Proposition} \label{fa28} Let $B_n$ act on a set $U${\rm .}
 Suppose we are given subsets $C_x\subset U$ {\rm (}$x\in \Omega${\rm ).}
\begin{itemize}
\ritem{(a)} If the inclusion $xC_y\subset C_{\LF(xy)}$ holds for all pairs
$(x,y)\in\Omega_1\times\Omega${\rm ,} then it
holds for all pairs in $B_n^+\times\Omega${\rm .}
\ritem{(b)} Assume the following\/{\rm :}
\begin{itemize} 
\ritem{(1)} The $C_x$ are nonempty and {\rm (}\/pairwise\/{\rm )} disjoint{\rm .}
\ritem{(2)} The properties of {\rm (a)} hold{\rm .}\end{itemize}
Then the $B_n$-action on $U$ is faithful{\rm .}\end{itemize}

\endproclaim

\demo{Proof} (a) We will show the desired result by induction on $\ell(x)$. If\break $\ell(x)\leq1$, there is nothing to
prove. Now let $\ell(x)>1$, say $x=uv$,\break $u,v\in B_n^+ -\{1\}$. Then $xC_y=u(vC_y)\subset
u(C_{\LF(vy)})\subset C_{\LF(u\,\LF(vy))}=C_{\LF(uvy)}=C_{\LF(xy)}$. (The two inclusions follow from the
induction hypothesis. The middle equality follows from (\ref{tf4}).) This proves the induction step and
thereby part (a).
\medbreak
(b) Let $\Sym(U)$ denote the group of permutations of $U$, and let $\pi\co B_n\ra \Sym(U)$ denote the action. Write $xu$ instead of $(\pi x)u$ ($x\in B_n$, $u\in U$). 
It is known that for any $z\in B_n$ there are $x,y\in B_n^+$ with $z=xy^{-1}$.
 Our proposition will therefore be proved if we show that for any $x,y\in B_n^+$, if $\pi(x)=\pi(y)$ then $x=y$. We will
show this by induction on $\ell(x)+\ell(y)$. 

Suppose $x,y\in B_n^+$ with $\pi(x)=\pi(y)$. If $\ell(x)+\ell(y)=0$ 
then $x=1$ and $y=1$, so certainly $x=y$. Consider now the case $\ell(x)+\ell(y)>0$. It is given that $C_1$ is
nonempty; choose any $u\in C_1$. By (a), we have $xu\in xC_1\subset C_{\LF(x)}$ and similarly $yu\in C_{\LF(y)}$. We
have $\pi x=\pi y$, whence $xu=yu$. It follows that $xu\in C_{\LF(x)}\cap C_{\LF(y)}$. By assumption (1), all
$C_z$ are disjoint however. It follows that $\LF(x)=\LF(y)$. Write $z=\LF(x)$, and define $x',y'\in B_n^+$ by
$x=zx'$, $y=zy'$. Note $z\neq 1$, because otherwise $x=y=1$, contradicting the fact that $\ell(x)+\ell(y)>0$.
It follows that
$\ell(x')+\ell(y')<\ell(x)+\ell(y)$. The induction assumption thus yields $x'=y'$ and hence $x=y$. This proves the
induction step and thereby part (b) of the proposition.\enddemo

The results in this section so far suffice to understand the
faithfulness proof in Sections~\ref{secd}, \ref{sece}. We now turn to some more combinatorial results which will be used in the proof of \ref{tf26}.

The {\it Charney length function\,} is the length function $\ell_\Omega\co B_n\ra\zz_{\geq0}$ with respect to
$\Omega$; i.e., $\ell_\Omega(x)$ is the smallest natural number $k$ such that there exist
$x_1,\ldots,x_k\in\Omega\cup\Omega^{-1}$ with $x=x_1\cdots x_k$. 

The center of $B_n$ is 
isomorphic to $\zz$ and, if $n\geq3$, generated by $\Delta^2$. We have a bijection $\zz\times(B_n^+ -\Delta
B_n^+)\ra B_n$ defined by $(k,x)\mapsto \Delta^kx$.

From the Artin presentation
 of the braid group, it follows that there exists an automorphism of $B_n$ which takes any $\sigma_i$ to its inverse.
We will denote this automorphism by $x\mapsto\barr{x}$.

The following theorem collects some combinatorial results.

\proclaimtitle{Garside, Thurston, Charney}
\proclaim{Theorem}  \label{fa74} \begin{itemize}
\ritem{(a)}  Let $(x_1,\ldots,x_k)$ denote the greedy form of some positive braid $x\in B_n^+${\rm .} Then
$\ell_\Omega(x)=k${\rm .}
\ritem {(b)}   Let $x\in B_n$. Then there are unique $y=y_x$ and $z=z_x$ both in $B_n^+$ with $x=y^{-1}z$
such that there is no $w\in \Omega_1$ with $\{y,z\}\subset wB_n^+${\rm .} They satisfy
$\ell_\Omega(x)=\ell_\Omega(y)+\ell_\Omega(z)${\rm .}

\ritem{(c)}   Let $x\in B_n^+ -\Delta B_n^+$ with $\ell_\Omega(x)=k${\rm .} Then $\ell_\Omega(\Delta^\ell
x)=\max(k+\ell,k,-\ell)$ for all $\ell\in\zz${\rm .}
\ritem{(d)}   Let $x\in B_n^+ -\Delta B_n^+$ with $\ell_\Omega(x)=k${\rm .} Then $\Delta^k \barr{x}\in B_n^+
-\Delta B_n^+$ and $\ell_\Omega(\Delta^k\barr{x})=k${\rm .}
\ritem{(e)}  The growth function
$$ \sum_{x\in B_n} z^{\ell_\Omega(x)}\in\zz[[z]] \label{fa81} $$
is rational.
\ritem{(f)}  There exists an algorithm that on input $n\in\zz_{\geq0}$ and $x\in B_n$ computes {\rm (}\/the
greedy forms of\/{\rm )} $y_x,z_x$ {\rm (}\/as defined in {\rm (b))} and $\ell_\Omega(x)${\rm .}
 The time the algorithm takes is
bounded by a polynomial in $n+\ell_\Omega(x)$.\end{itemize}

\endproclaim

  Charney's result Theorem~\ref{fa74}(e) becomes even more remarkable if one knows that for most other finite generating
subsets of $B_n$ (including the Artin generating set $\{\sigma_1,\ldots,\sigma_{n-1}\}=\Omega_1$) it is still
unknown whether the growth function with respect to it is rational. This should not be confused with Deligne's result
\cite{del} that the growth function of positive braids
$$ \sum_{x\in B_n^+}z^{\ell(x)} $$
is rational (see also \cite{xu}). 

In contrast to Theorem~\ref{fa74}(f) Paterson and Razborov \cite{pr} proved that computing the length of a braid in $B_n$ with
respect to $\Omega_1$ (with $n$ variable) is an NP-complete problem.

\section{The representation} \label{secc}
\advance\eqcount by 3

Let $R$ denote a commutative ring and $q,t\in R$ two invertible
 elements. Let $V$ denote the free $R$-module with basis $\{x_s\mid s\in\Ref\}$. Thus, the dimension of $V$ is
$m:=|\Ref|=n(n-1)/2$. We will also write $x_{ij}$ instead of $x_{s(i,j)}$ where $1\leq i<j\leq n$. We define a
representation $\rho\co B_n\ra \GL(V)$ as follows (action of $\GL(V)$ on $V$ from the left; instead of $(\rho x)v$,
we use the simpler notation $xv$, $x\in B_n$, $v\in V$):
$$ \begin{array}{ll}
\sigma_k\,x_{k,k+1} =tq^2\,x_{k,k+1};   \\[4pt]
\sigma_k\,x_{ik}=(1-q)\,x_{ik}+q\,x_{i,k+1}, & i<k;   \\[4pt]
\sigma_k\,x_{i,k+1}=x_{ik}+ tq^{k-i+1}(q-1)\,x_{k,k+1},\ \ & i<k;   \\[4pt]
\sigma_k\,x_{kj}=tq(q-1)\,x_{k,k+1}+q\,x_{k+1,j}, & k+1<j;   \\[4pt]
\sigma_k\,x_{k+1,j}=x_{kj}+(1-q)\,x_{k+1,j}, & k+1<j;   \\[4pt]
\sigma_k\,x_{ij}=x_{ij}, & i<j<k\mbox{ \ or \ }k+1<i<j;  \\[4pt]
\sigma_k\,x_{ij}=x_{ij}+tq^{k-i}(q-1)^2\,x_{k,k+1}, & i<k<k+1<j. \end{array} $$
It should be proved here that these formulas do indeed define a representation, i.e., that they respect relations
(\ref{tf1}) and (\ref{tf2}) in the Artin presentation of the braid group, and that $\rho\sigma_k$ is invertible. This is a
straightforward though tedious task which we leave to the reader. 

\demo{{R}emark} In \cite{kra} the author uses another basis $\{v_{ij}\mid 1\leq i<j\leq n\}$ of the same module $V$.
Its relation with $\{x_{ij}\}_{ij}$ is given by
\be v_{ij}=x_{ij}+(1-q)\sum_{i<k<j}x_{kj},\ \ \ \ 
   x_{ij}=v_{ij}+(q-1)\sum_{i<k<j}q^{k-1-i}\,v_{kj}. \label{tf25} \ee
In \cite{kra} one can also find a topological interpretation of $v_{ij}$. 
Combination with (\ref{tf25}) then results in a topological interpretation of $x_{ij}$. In the present paper, we will
not consider any other bases than $\{x_{ij}\}_{ij}$ and its dual. A quicker proof of our formulas defining a
representation is obtained if one is willing to accept the formulas with respect to $\{v_{ij}\}$ in \cite{kra}, by
combining with (\ref{tf25}).\enddemo

Let $V^*$ denote the dual of $V$ and let 
$\langle \cdot,\cdot\rangle \co V^*\times V\ra R$ denote the pairing. Let $\{x^*_s\mid s\in\Ref\}$ denote the basis
dual to $\{x_s\}_s$, and write $x^*_{ij}=x^*_{s(i,j)}$. In formula:
$$ V^*=\bigoplus_{s\in\Ref}Rx^*_s,\ \ \ \ \langle x^*_r,x_s\rangle =\delta_{rs}. $$
Let $\End(V)$ act on $V^*$ on the right by $\langle uA,v\rangle =\langle u,Av\rangle $ for all $A\in \End(V)$, $(u,v)\in V^*\times V$. Again, we will write $vx$ instead of $v(\rho x)$ ($x\in B_n$, $v\in V^*$).

The braid group action on $V^*$ is given by
\begin{eqnarray}\quad
 x^*_{k,k+1}\sigma_k&\hskip-9pt = \hskip-9pt& t\Bigg[q^2\,x^*_{k,k+1}\ +\ q(q-1)\sum_{k+1<b}x^*_{kb}
 \label{tf35} \\ && +  (q-1)\sum_{a<k}q^{k-a+1}\,x^*_{a,k+1} \ +\ (q-1)^2\sum_{a<k<k+1<b}q^{k-a}x^*_{ab}
\Bigg] \nonumber
\end{eqnarray}
and
\be x^*_{ij}\sigma_k=  \left\{ \begin{array}{ll}
(1-q)\,x^*_{ik}+x^*_{i,k+1}, & i<k,\ j=k;   \\[4pt]
q\,x^*_{ik}, & i<k,\ j=k+1;  \\[4pt]
x^*_{k+1,j}, & i=k,\ j>k+1;   \\[4pt]
(1-q)\,x^*_{k+1,j}+q\,x^*_{kj}, & i=k+1,\ j>k+1;   \\[4pt]
x^*_{ij}, & \{i,j\}\cap\{k,k+1\}=\emptyset.\ \end{array} \right. \label{tf36} \ee

The results in this section so far are sufficient background for
reading the faithfulness proof in the next two sections. The remainder
of this section deals with a few identities which will be used in Section~\ref{secf} to prove some more properties of the representation.

We define a linear map $T(u)\co V\ra V$ depending on a parameter $u$ by
\begin{eqnarray*}
T(u)\co x_{ij}&\mapsto&
\sum_{i<k<\ell<j}(1-u)^2u^{i+\ell}\, x_{k\ell}\ \ +
\sum_{i=k<\ell<j}(1-u)u^{i+\ell}\, x_{k\ell}\\
&&+\  \sum_{i<k<\ell=j}(1-u)u^{i+\ell}\, x_{k\ell}\ \ +
\sum_{i=k<\ell=j}u^{i+\ell}\, x_{k\ell},
\end{eqnarray*}
each sum ranging over those $(k,\ell)$ with $1\leq k<\ell\leq n$ satisfying the indicated inequalities and identities.

An obvious total ordering on the basis 
$\{x_{ij}\}$ makes $T(u)$ into a triangle matrix with powers of $u$ on the diagonal. In particular, $T(u)$ is invertible
if $u$ is. (A fact which we will not need is that $T(u)T(u^{-1})=1$.)

In order to indicate the dependence 
of $\rho$ on $q,t$, we write $\rho(x,q,t)$ ($x\in B_n$). Thus, for any two invertible elements $q',t'$ in some ring,
$\rho(x,q',t')$ is obtained from the matrix of $\rho x$ with respect to $\{x_{ij}\}_{ij}$ by entry-wise replacing $q$ by
$q'$ and $t$ by $t'$.

Recall the automorphism of $B_n$ denoted $x\mapsto\barr{x}$ and mapping each $\sigma_i$ to its inverse.

\proclaim{Lemma} \label{tf8} 
For all $x\in B_n${\rm ,} one has \, $T(q)\,\rho(x,q^{-1},t^{-1})\,T(q)^{-1}=\rho(\barr{x},q,t)${\rm .}\endproclaim

\demo{Proof} The proof is straightforward and left to the reader. We only remark that it suffices to give a proof for
$x\in\Omega_1$, because both sides of the desired identity are group homomorphism images of $x$.\enddemo

If the braid group is viewed as the mapping class group of
 the punctured disk, all punctures being real, then the matrix $T$ corresponds to complex conjugation. For the
application of \ref{tf8} we have in mind (\ref{tf26}), it is important that $T(q)$ does not involve~$t$.

\proclaim{Lemma} \label{tf24} 
We have \ $\Delta\,x_{n+1-j,n+1-i}=tq^{i+j-1}\,x_{ij}$ \ whenever $1\leq i<j\leq n${\rm .}\endproclaim

\demo{Proof} Perhaps the best way to prove it is by having a topological interpretation of the representation. (See for example \cite{law} or \cite{kra}.) We will follow a more elementary path instead, which completely avoids topological arguments. Basically the idea is to multiply the matrices $\rho\sigma_i$ according to a factorization of $\Delta$ like $\Delta=(\sigma_1\cdots\sigma_n)(\sigma_1\cdots\sigma_{n-1})\cdots(\sigma_1\sigma_2)\sigma_1$. We reduce the amount of calculations as follows.

Define $A\in\GL(V)$ by the expected formula for $\Delta$:  $A\,x_{n+1-j,n+1-i}=tq^{i+j-1}\,x_{ij}$. Our goal is then to prove $A=\rho\Delta$.

\medbreak {\it Claim} 1.  $A(\rho\sigma_k)A^{-1}=\rho\sigma_{n-k}$ whenever $1\leq k\leq n-1$. This is a
straightforward computation, which will be left to the reader.

\medbreak {\it Claim} 2.  The centralizer in $\GL(V)$ of $\rho B_n$ consists of the scalar matrices only. (This is
closely related to the irreducibility of $V$, which was established by Zinno \cite{zin}.) In order to prove Claim~2, let
$B\in\GL(V)$ commute with each element of $\rho B_n$. We must then show $B$ to be a scalar matrix. Note that
$x_{k,k+1}$ is an eigenvector of $\sigma_k$ with eigenvalue $tq^2$. From our formula for the $\sigma_k$-action on
$V^*$, it readily follows that the eigenvalue $tq^2$ is simple; indeed, all remaining eigenvalues depend on $q$ only.
Since $v_{12}$ is an eigenvector of $\sigma_1$ with simple eigenvalue, and $B$ commutes with $\sigma_1$, one has
$Bx_{12}=\lambda x_{12}$ for some invertible $\lambda\in R$. After multiplying $B$ with a scalar matrix, we may
assume
$Bx_{12}=x_{12}$; our task is then to show $B=1$. Using the identities
$\sigma_k\,x_{1k}=(1-q)\,x_{1k}+q\,x_{1,k+1}$ ($1<k$) one inductively finds $Bx_{1j}=x_{1j}$. Using the identities
$\sigma_k\,x_{k+1,j}=x_{kj}+(1-q)\,x_{k+1,j}$ ($k+1<j$), one inductively finds $Bx_{ij}=x_{ij}$.
 This shows that $B=1$, thus proving Claim~2.

Note $\Delta\sigma_k\Delta^{-1}=\sigma_{n-k}$ for all $k$. Thus, the property of $A$ formulated in Claim~1 is also
satisfied by $\rho\Delta$. In other words, $A^{-1}(\rho\Delta)$ is in the centralizer of $\rho B_n$. By Claim~2,
$A^{-1}(\rho\Delta)$ is scalar, so it suffices to show that $Ax_{n-1,n}=\Delta x_{n-1,n}$. We use the following
factorization of $\Delta$:
$$ \Delta=\sigma_1(\sigma_2\sigma_1)(\sigma_3\sigma_2)\cdots(\sigma_{n-1} \sigma_{n-2}) \Delta_{n-2} $$
where 
$$ \Delta_{n-2}=(\sigma_1\cdots \sigma_{n-2})(\sigma_1\cdots \sigma_{n-3})\cdots (\sigma_1\sigma_2)\sigma_1. 
$$ First of all, $\Delta_{n-2}x_{n-1,n}=x_{n-1,n}$. Moreover, for $3\leq k\leq n$, one has
\begin{eqnarray*} \sigma_{k-1}\sigma_{k-2}x_{k-1,k}&=&\sigma_{k-1}(x_{k-2,k}+(1-q)\,x_{k-1,k})\\
 & =&(x_{k-2,k-1}+tq^2(q-1)\,x_{k-1,k})\\
&&+\ (1-q)tq^2x_{k-1,k}=x_{k-2,k-1}. \end{eqnarray*}
It follows that $\Delta x_{n-1,n}= \sigma_1(\sigma_2\sigma_1)(\sigma_3\sigma_2)\cdots(\sigma_{n-1} \sigma_{n-2})\Delta_{n-2}x_{n-1,n}=
\sigma_1x_{12}=tq^2x_{12}=Ax_{n-1,n}$. This finishes the proof.\enddemo

\section{Faithfulness} \label{secd}
\advance\eqcount by 6

Recall our representation $\rho\co B_n\ra\GL(V)$. We often tacitly identify $\rho x$
 with its matrix with respect to the basis $\{x_{ij}\}_{ij}$. 

Observe that for any positive braid $x$, the entries of the matrix of
$\rho x$ are in $\zz[q,q^{-1},t]$. This follows from the matrices given
in Section \ref{secc} and the fact that $B_n^+$ is generated by $\Omega_1=\{\sigma_1,\ldots,\sigma_{n-1}\}$.

From now on, we take the one-variable Laurent polynomial ring $R=\rr[t^{\pm1}]$ for base ring, and $q\in\rr\subset R$ with $0<q< 1$. Put
$$ V_1:=\bigoplus_{s\in\Ref}\ \rr[t]\,x_s\subset V=\bigoplus_{s\in\Ref}\ \rr[t^{\pm1}]\,x_s.$$
We thus have $B_n^+ V_1\subset V_1$.

Note that all entries of $\rho\sigma_k$ are in $\{0,1,q,1-q\}+t\,\zz[q,q^{-1},t]$. By our assumption that $q$ is a real number with $0<q<1$, they are in $\rr_{\geq0}+t\,\rr[t]$. On putting
\begin{eqnarray*} V_2&=&\bigoplus_{s\in\Ref}\Bigl( \rr_{\geq0}+t\,\rr[t] \Bigr) \,x_s=\Biggl(
\bigoplus_{s\in\Ref}\rr_{\geq0}\,x_s\Biggr) \oplus\,tV_1\\
& =&\Big\{ v\in V\ {\Big\vert}\ \mbox{for all }s\in\Ref\co
\langle x_s^*,v\rangle \in\rr_{\geq0}+t\,\rr[t]\Big\}, \end{eqnarray*}
 we have $B_n^+V_2\subset V_2$. 

\demo{Definition} For $A\subset\Ref$, define
$$ D_A=\Big\{ v\in V_2\ {\Big\vert}\ \mbox{for all }s\in\Ref\co \langle x_s^*,v
\rangle \in t\,\rr[t]\Leftrightarrow s\in A\Big\}. $$

\hglue18pt Thus, $V_2$ is the disjoint union of the $D_A$ ($A\subset\Ref$).
\enddemo

Let $x\in B_n^+$, $A\subset\Ref$. Then there is a unique $B\subset\Ref$ with $xD_A\subset D_B$. (A formula
for $B$ is given in \ref{tf17}.) Notation: $B=xA$. Let $2^\Ref$ denote the power set of $\Ref$. The map
$B_n^+\times 2^\Ref\ra2^\Ref$, $(x,A)\mapsto xA$ defines an action of $B_n^+$ on $2^\Ref$. (This follows
from the facts that $\rho$ is a representation and that $\rho B_n^+$ preserves $V_2$.) An explicit formula
for the $B_n^+$-action on $2^\Ref$ is as follows.

\proclaim{Lemma} \label{tf17} Let $A\subset \Ref$ and $1\leq k\leq n-1${\rm .}
 Then $\sigma_kA$ is the set of those $s(i,j)$ with $1\leq i<j\leq n$ and
$$ \left\{ \begin{array}{ll}
\mbox{\rm true statement,} & i=k,\ j=k+1; \\
\{s(i,k),s(i,k+1)\}\subset A, & i<k,\ j=k; \\
s(i,k)\in A, & i<k,\ j=k+1; \\
s(k+1,j)\in A, & i=k,\ j>k+1; \\
\{s(k+1,j),s(k,j)\}\subset A,\ \ & i=k+1,\ j>k+1; \\
s(i,j)\in A, & \{i,j\}\cap\{k,k+1\}=\emptyset.\end{array} \right.$$\endproclaim

\demo{Proof} This is readily obtained from the formulas 
for the $\sigma_k$-action on $V^*$, (\ref{tf35}) and (\ref{tf36}).\enddemo

If one were given the formula of Lemma~\ref{tf17} only, it would not be obvious that it defines a $B_n^+$-action on $2^{\Ref}$; but we get the proof of this fact for free as a consequence of our representation.

The notation $xA$ ($A\subset\Ref$) can have rather different meanings according to whether $x\in S_n$ or $x\in B_n^+$. If $x\in S_n$ then $xA=\{xa\mid a\in A\}$, involving multiplication in the symmetric group. For $x\in B_n^+$, the notation refers to the $B_n^+$-action on $2^\Ref$.

Note that the $B_n^+$-action on $2^\Ref$ preserves inclusions, i.e., for $A\subset B\subset\Ref$ and $x\in B_n^+$ one has $xA\subset xB$.

\demo{Definition}  We define a map $\Inv\co S_n\ra2^\Ref$ by
$$ \Inv(x)=\Big\{s(i,j)\ {\Big\vert}\ 1\leq i<j\leq n,\ x^{-1}i>x^{-1}j\Big\}. $$

Note:  $\ell(x)=|\Inv(x)|$. Moreover, for all $x,y\in S_n$, we have 
\begin{eqnarray}
 x\leq xy&\Longleftrightarrow& \ell(xy)=\ell(x)+\ell(y)\label{tf34} \\
&\Longleftrightarrow &
\Inv(xy)=\Inv(x)\cup x\,
\Inv(y)\,x^{-1}\Longleftrightarrow \Inv(x)\subset \Inv(xy).   \nonumber\end{eqnarray}

The image of $\Inv$ will be denoted by $\Inv(S_n)$. As $\Inv$ is injective, one may identify $S_n$ with $\Inv(S_n)$. We will however distinguish them in our notation, because otherwise it would cause confusion.

It can be shown that the $B_n^+$-action on $2^\Ref$ does not preserve $L(S_n)$. (See \ref{tf45} for more details.) In particular, the obvious definition $C_x=D_{Lr^{-1}x}$(?) does not satisfy the conditions of \ref{fa28}(b). Therefore we need a new idea, which is as follows.
\enddemo

\demo{Definition} A set $A\subset\Ref$ is said to be a {\it half-permutation\,}\footnote{Half-permutations
are more commonly called {\it closed sets}.} if, whenever
$1\leq i<j<k\leq n$, one has
$$ s(i,j),\ s(j,k)\in A\Longrightarrow s(i,k)\in A. $$
We will denote the set of half-permutations by $\HP$.

Every element of $L(S_n)$ is a half-permutation.
 A fact which we shall not need is that a subset $A\subset\Ref$ is in $\Inv(S_n)$ if and only if both $A$ and
$(\Ref-A)$ are half-permutations. This explains the terminology of half-permutations. There is another interpretation
of half-permutations as follows. There is a bijection from the set of (partial) orderings $<_0$ on $I_n=\{1,\ldots,n\}$
with $(i<_0j)\Rightarrow\break (i<j)$, to $\HP$, which takes $<_0$ to $\{s(i,j)\mid i<_0j\}$.
\enddemo

We next record a few combinatorial results whose proofs are deferred to the next section for the sake of readability.

\proclaim{Lemma}  \label{tf9} Let $x\in B_n^+${\rm ,} $A\in\HP$. Then $xA\in \HP${\rm .}\endproclaim

\demo{Proof} For a proof, see \ref{tf13}.\enddemo

Recall that a {\it greatest\,} element in a (partially) ordered set is an element greater than all other elements.

\proclaim{Lemma/Definition}  \label{tf10} For every half\/{\rm -}\/permutation $A$ there is a greatest {\rm
(}\/with respect to inclusion\/{\rm )}
$B\in\Inv(S_n)$ with $B\subset A${\rm .} Notation\/{\rm :} $B=\Pro(A)${\rm .}\endproclaim

\demo{Proof} See \ref{tf14}.\enddemo

\demo{Definition}  Let $\,\GB$ (Greatest Braid) denote the map 
$$ \GB=r\,\Inv^{-1}\,\Pro\co \HP\ra\Omega. $$
Moreover, for $x\in\Omega$, define 
$$ C_x=\cup\,\Big\{D_A\ {\Big\vert}\ A\in \HP,\ \GB(A)=x \Big\}. $$

Recall the $B_n^+$-action on $\Omega$ defined by $B_n^+\times\Omega\ra\Omega$, $(x,y)\mapsto\LF(xy)$.
\enddemo

\proclaim{Lemma} \label{tf11} The map $\,\GB\co \HP\ra\Omega$ is $B_n^+$\/{\rm -}\/equivariant{\rm .}
 In formula{\rm ,} if $x\in B_n^+${\rm ,} $A\in\HP${\rm ,} $y=\GB(A)${\rm ,} then $\GB(xA)=\LF(xy)${\rm .}\endproclaim

\demo{Proof} See \ref{tf15}.\enddemo

Since $\GB\co \HP\ra\Omega$ is surjective (indeed, $\GB(\Inv(r^{-1}x))=x$ for all $x\in\Omega$), one may call $\HP$ a refinement of $\Omega$.

\proclaim{Lemma} \label{tf12} \begin{itemize}
\ritem{(a)} The $C_x$ are disjoint and nonempty{\rm .}
\ritem{(b)} We have $xC_y\subset C_{\LF(xy)}$ for all $(x,y)\in B_n^+\times\Omega${\rm .}\end{itemize}\endproclaim

{\it Proof}. (a) The set $C_x$ is nonempty 
because $\emptyset\neq D_{\Inv(r^{-1}x)}\subset C_x$. That the $C_x$ are disjoint is trivial.

(b) The definition of $C_y$ reads $C_y=\cup\,\{D_A\mid A\in \HP,\ \GB(A)=y\}$. 
We must therefore show $xD_A\subset C_{\LF(xy)}$ whenever $A\in\HP$, $\GB(A)=y$. Our proof is summarized by the
following chain:
$$ xD_A\stackrel{(1)}{\subset} D_{xA}\stackrel{(2)}{\subset}
   C_{\GB(xA)}\stackrel{(3)}{=} C_{\LF(xy)}. $$
Here, (1) follows from the definition of the $B_n^+$-action on $2^\Ref$. 
In order to justify (2), note that $xA\in \HP$ by Lemma \ref{tf9}. By \ref{tf10} then, $\Pro(xA)$ is defined and hence
so are $\GB(xA)$ and the right-hand side of (2). Inclusion (2) follows by definition of $C_z$. Identity (3) is
Lemma~\ref{tf11}. This finishes the proof of (b).\hfill\qed

 \proclaim{Theorem} \label{tf7} The representation $\rho\co B_n\ra\GL(V)$ is faithful{\rm .}\endproclaim

{\it Proof}. This is an immediate consequence of \ref{tf12} and \ref{fa28}(b) (with $U=~V$).
\phantom{almost done}\hfill\qed
\medbreak
The considerations of this section can be illustrated by the commutative diagram of $B_n^+$-equivariant maps in
Figure~1 below. The arrows pointing to the left are inclusions.

\begin{center}
\BoxedEPSF{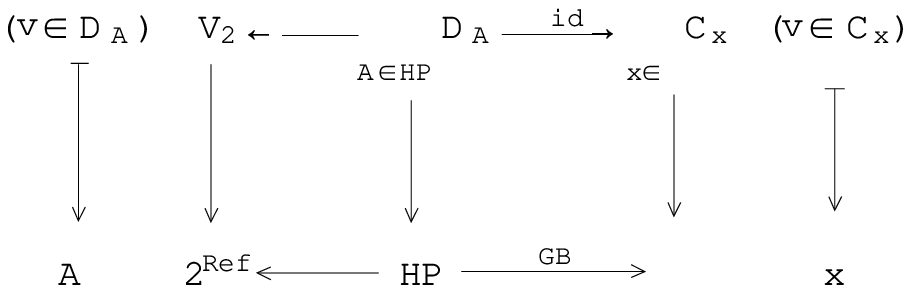 scaled1000}
\end{center}
\vglue6pt
\centerline{Figure 1. A $B_n^+$-equivariant commutative diagram}

\numbereddemo{{R}emark} \label{tf45} The (easy) proof of \ref{tf12}
 shows a more general statement as follows. Let $\HP'\subset2^\Ref$ denote a $B_n^+$-invariant subset, and $\GB'\co
\HP'\ra\Omega$ a surjective $B_n^+$-equivariant map. Then the sets $C_x':=\mathbold{\cup}\{D_A\mid A\in\HP',\
\GB'(A)=x\}$ satisfy the same conclusion of \ref{tf12}, thus proving once more that $\rho$ is faithful. In this section,
we have constructed one such a pair $(\HP',\GB')$, namely, $(\HP,\GB)$. The following questions arise. Are there more
such pairs $(\HP',\GB')$? Is there a best pair, whatever that means? Any pair with the desired properties cannot
involve
$\HP'=\Inv(S_n)$, because $\Inv(S_n)$ is not $B_n^+$-invariant, as is proved by the following counterexample: If
$n=4$, then $A:=\{s_{13},s_{14},s_{23},s_{24}\}\in\Inv(S_n)$, but $\sigma_2A=\{s_{23},s_{14}\}\not\in\Inv(S_n)$.
Another solution is given by $\HP'=\HP_0$, the smallest $B_n^+$-invariant subset of $\HP$ containing $\Inv(S_n)$,
and $\GB'=\GB_0$, the appropriate restriction of $\GB$. For example, if $n=4$, then $|\HP_0|=25$ and
$\HP_0=\Inv(S_n)\cup \{\{s_{23},s_{14}\}\}$. The sets $\HP_0$ seem to be rather messy, and $(\HP,\GB)$ is a
comfortable solution after all.\enddemo

 \vglue-8pt
\section{Half-permutations} \label{sece}
\advance\eqcount by 7
\vglue-4pt

The aim of this section is to prove some combinatorial results among which are the promised lemmas of the previous section.

\numbereddemo{{R}emark} \label{tf16}
 There exists an involutory automorphism of the system $(S_n,\ S,\ B_n^+,$ $r\co S_n\ra B_n^+,\ \Ref,\ B_n^+\times
2^\Ref\ra2^\Ref,\ \HP,\ \Inv,\  \Pro,\ \GB)$, defined by conjugation by $w_0$ in $S_n$, or by $\Delta$ in $B_n$. (Note:
$\Delta^2$ is central in $B_n$.) The easy proof is left to the reader. For example, this involution maps $s_k$ to
$s_{n-k}$. Especially the fact that the involution preserves the $B_n^+$-action on $2^\Ref$ is remarkable. This
symmetry will prove useful as it can be used to reduce the number of cases in a few case-by-case proofs.
\enddemo

 \proclaim{Lemma} \label{tf13} Let $x\in B_n^+$, $A\in\HP$. Then $xA\in \HP$.\endproclaim

{\it Proof}. One may suppose $x\in\Omega_1$, say $x=\sigma_k$. Let $1\leq p<q<r\leq n$. We must prove (H):
$s(p,q),\ s(q,r)\in xA\Rightarrow s(p,r)\in xA$. Modulo the symmetry of \ref{tf16}, there are five cases to
consider, as shown in the first two columns of Figure~2. For each of these cases, the table in
Figure~2 gives a statement in terms of $A$ which is equivalent to a given one among $s_{pq}\in
xA$, $s_{qr}\in xA$, $s_{pr}\in xA$. This table is a consequence of \ref{tf17}. Using the table, one readily
verifies (H). As an example, we do Case 4: 
\begin{eqnarray*}
 s(p,q),\ s(q,r)\in xA &\Rightarrow& s(p,q),\ s(p,k+1),\
 s(k+1,r)\in
A \\
&\Rightarrow & s(p,k+1),\
 s(k+1,r)\in A \\
&\Rightarrow& s(p,r)\in A\, \Rightarrow\, s(p,r)\in xA ,
\end{eqnarray*} which proves
case 4.\hfill\qed

\vspace*{10pt} \begin{center}
\makebox{\begin{tabular}{|c|l||l|l|l|} \hline
& \ruup & \parbox[t]{27mm}{$s(p,q)\in xA\\ \Longleftrightarrow \ldots$\rudown} &
\parbox[t]{27mm}{$s(q,r)\in xA\\ \Longleftrightarrow \ldots$} &
\parbox[t]{27mm}{$s(p,r)\in xA\\ \Longleftrightarrow \ldots$} \\ \hline\hline
1\ruup & \parbox[t]{27mm}{$\{p,q,r\}\cap\\ \{k,k+1\}=\emptyset$\rudown} &
$s(p,q)\in A$ & $s(q,r)\in A$ & $s(p,r)\in A$ \\ \hline
2\ruup & $r=k$ & $s(p,q)\in A$ & \parbox[t]{29mm}{$s(q,r)\in A$ \ and\\ $s(q,k+1)\in A$\rudown} & \parbox[t]{29mm}{$s(p,r)\in A$ \ and\\ $s(p,k+1)\in A$} \\ \hline
3\ruboth & $q<k,\ r=k+1$ & $s(p,q)\in A$ & $s(q,k)\in A$ & $s(p,k)\in A$ \\ \hline
4\ruup & $q=k,\ r>k+1$ & \parbox[t]{29mm}{$s(p,q)\in A$ \ and\\ $s(p,k+1)\in A$\rudown} & $s(k+1,r)\in A$ & $s(p,r)\in A$ \\ \hline
5\ruup & $q=k,\ r=k+1$ & \parbox[t]{29mm}{$s(p,q)\in A$ \ and\\ $s(p,r)\in A$\rudown} & true & $s(p,q)\in A$ \\ \hline
\end{tabular}} \end{center} 
\centerline{Figure 2. To the proof of \ref{tf13} \label{tffig1}}

\proclaim{Lemma} \label{tf18} Let $A\in\HP$, $x\in S_n$, $\Inv(x)\subset A$, $B=x^{-1}(A-\Inv(x))x$. Then $B\in \HP$.\endproclaim

{\it Proof}. First we prove the lemma for the case $\ell(x)=1$, say $x=s_k$. Notice that $x^2=1$. Let $1\leq
p<q<r\leq n$. We must show (H): $s_{pq},\ s_{qr}\in B\Rightarrow s_{pr}\in B$. First consider the case where
$\{p,q,r\}\cap\{k,k+1\}$ consists of at most one element. Write $(p',q',r')=(xp,xq,xr)$. Then $p'<q'<r'$. From
$s(p,q)\in B$ we find $s(p',q')=x\,s(p,q)\,x\in A$; similarly $s(q,r)\in B$ implies $s(q',r')=x\,s(q,r)\,x\in A$. As
$A$ is a half-permutation and $p'<q'<r'$, it follows that $s(p',r')\in A$. Hence $s(p,r)=x\,s(p',r')\,x\in B$, thus
proving (H) if $|\{p,q,r\}\cap\{k,k+1\}|\leq 1$. Because of the symmetry of \ref{tf16}, it remains only to
consider the case $q=k$, $r=k+1$. Then the left-hand side of (H) implies $s_{qr}\in B$, whence $x\in B$,
whence $x\in xBx=A-\Inv(x)=A-\{x\}$, a contradiction. This proves (H) in the case $q=k$, $r=k+1$, thus
establishing the lemma for $\ell(x)=1$.

We finish the proof of the lemma by induction on $\ell(x)$. For $\ell(x)\leq1$ there is nothing left to prove. Suppose $u\leq uv=x$ with $u,v\in S_n-\{1\}$. Recall (\ref{tf34}) that $\Inv(x)$ is the disjoint union of $\Inv(u)$ with $u\,\Inv(v)\,u^{-1}$. Since $\Inv(x)\subset A$, we have $\Inv(u)\subset A$. Applying the induction hypothesis to $(A,u)$ shows that $C:=u^{-1}(A-\Inv(u))u$ is a half-permutation. From $\Inv(x)\subset A$ we find $\Inv(v)\subset C$. Applying the induction hypothesis to $(C,v)$ then yields
\begin{eqnarray*}
 \HP\ni v^{-1}(C-\Inv(v))v&=&v^{-1}(u^{-1}(A-\Inv(u))u-\Inv(v))v\\
&
=&v^{-1}u^{-1}(A-\Inv(u)-u\,\Inv(v)\,u^{-1})uv\\
&=&x^{-1}(A-\Inv(x))x=B. 
\end{eqnarray*} This proves the induction step and
hence the lemma.\hfill\qed

\proclaim{Lemma/Definition} \label{tf14} For every half\/{\rm -}\/permutation $A$ there is a greatest {\rm (}\/with
respect to inclusion\/{\rm )}
$B\in\Inv(S_n)$ with $B\subset A${\rm .} Notation\/{\rm :} $B=\Pro(A)${\rm .}\endproclaim

{\it Proof}.
Recall (\ref{tf34}) that for $x,y\in S_n$ we have $x\leq y\Leftrightarrow\Inv(x)\subset \Inv(y)$. So an equivalent
formulation of the lemma is that $P:=\{y\in S_n\mid \Inv(y)\subset A\}$ contains a greatest element. This is the
formulation which we will prove.

Note that the ordering on $P$ is generated by $x\leq xs$ whenever true, with $x,xs\in P$, $s\in S$. Let $x,xs,xt\in P$ with $x\leq xs$, $x\leq xt$, $s,t\in S$ ($s\neq t$). Since $P$ is finite and has a smallest element, it suffices (by a well-known elementary result on partial ordered sets) to show that there exists then $y\in P$ with $xs,xt\leq y$. Let $m_{st}\in\{2,3\}$ denote the order of $st$, and put $y=xst$ if $m_{st}=2$, and $y=xsts$ if $m_{st}=3$. It is well-known that $xs,xt\leq y$. We claim that $y\in P$. 
The lemma would clearly follow from this claim. We consider two cases according to the value of $m_{st}$.

\vglue5pt {\it Case} {\rm 1}. $m_{st}=2$. Then $\Inv(y)=\Inv(xs)\cup\Inv(xt)\subset A$ whence $y\in
P$. 

\vglue5pt {\it Case}  {\rm 2}. $m_{st}=3$. Write $s=s_k$, $t=s_{k+1}$. Define $C=x^{-1}(A-\Inv(x))x$. Let
$\amalg$ denote disjoint union. For any $u\in S_n$ with $x\leq xu$ we have $xu\in P\Leftrightarrow
\Inv(xu)\subset A\Leftrightarrow
\Inv(x)\amalg x\,\Inv(u)\,x^{-1}\subset A \Leftrightarrow x\,\Inv(u)\,x^{-1}\subset A-\Inv(x)\Leftrightarrow
\Inv(u)\subset x^{-1}(A-\Inv(x))x=C$, thus showing (for any $u\in S_n$ with $x\leq xu$):
\be xu\in P\Longleftrightarrow \Inv(u)\subset C.\label{tf19} \ee
Applying (\ref{tf19}) to $u=s,t$ gives $s,t\in C$; i.e., $s(k,k+1),\ s(k+1,k+2)\in C$. By \ref{tf18}, we have $C\in\HP$,
which means that we may conclude $s(k,k+2)\in C$. Hence $\Inv(sts)=\{s(k,k+1), s(k+1,k+2),s(k,k+2)\}\subset C$.
Applying (\ref{tf19}) in the reverse direction to $u=sts$ we find $y=xsts\in P$. This finishes case 2 and
thereby the proof of the lemma.\hfill\qed \medbreak

Notice that $\,\Pro\,$ is a projection; i.e., $\Pro^2=\Pro$. Moreover, $\Pro\,\Inv=\Inv$.\pagebreak

As to the following lemma, we will only make use of the special case of (a) where $\ell(x)=1$. We prove the entire lemma because it appears to have interest of its own. Recall the $B_n^+$-action on $2^\Ref$ defined in the previous section (or by \ref{tf17}), and which preserves $\HP$ by \ref{tf13}.

\proclaim{Lemma} \label{tf20}
\begin{itemize}
\ritem{(a)} Let $A\in\HP$, $x\in S_n${\rm .} Then $(rx)A$ equals the greatest {\rm (}\/with respect to inclusion\/{\rm
)} half\/{\rm -}\/permutation
$B$ with 
\be \Inv(x)\subset B\subset\Inv(x)\cup xAx^{-1}. \label{tf21} \ee
{\rm (}\/In particular{\rm ,} a greatest such half\/{\rm -}\/permutation exists{\rm .)}
\ritem{(b)} For $x,y\in S_n$ with $x\leq xy$ we have $(rx)\Inv(y)=\Inv(xy)${\rm .} In particular {\rm (}\/for
$y=1$\/{\rm ),}
$(rx)\emptyset=\Inv(x)${\rm .}\end{itemize}

\endproclaim

\demo{Proof} We start by proving (a) if $\ell(x)=1$. Write $x=s_k$, and note $\Inv(x)=\{x\}$. By \ref{tf13}, we have $(rx)A\in \HP$. From the definition of $(rx)A$ one readily finds $\{x\}\subset (rx)A\subset\{x\}\cup xAx$.

It remains to show, for any half-permutation $B$, that (\ref{tf21}) implies
 $B\subset (rx)A$. Suppose $s_{ij}\in B$, $1\leq i<j\leq n$. We must prove $s_{ij}\in (rx)A$. We consider four cases.

\medbreak{\it Case} 1.  $i=k$, $j=k+1$. Then $s_{ij}\in(rx)A$ by \ref{tf17}.

\medbreak{\it Case} 2.  $i<k$, $j=k+1$. We have $x\neq s_{ij}\in B\subset\{x\}\cup xAx$, whence $s_{ij}\in xAx$,
whence $s_{ik}=xs_{ij}x\in A$, whence $s_{ij}=s_{i,k+1}\in (rx)A$ by \ref{tf17}.

\medbreak{\it Case} 3. $i<k$, $j=k$. Then similarly to Case 2, we have $x\neq s_{ij}\in B\subset\{x\}\cup
xAx$, whence $s_{ij}\in xAx$ and $s_{i,k+1}=xs_{ij}x\in A$. Moreover, as $s_{ik},\ s_{k,k+1}\in B$ and
$B\in\HP$, we also have $s_{i,k+1}\in B$. In Case 2 we already saw that $s_{i,k+1}\in B$ implies $s_{ik}\in
A$. Summarizing, we have
$s_{ik},s_{i,k+1}\in A$ whence $s_{ij}=s_{ik}\in(rx)A$ by \ref{tf17}. This finishes Case 3.

\medbreak{\it Case} 4. $\{i,j\}\cap\{k,k+1\}=\emptyset$. Then $s_{ij}\in B$ readily implies $s_{ij}\in A$ and hence
$s_{ij}\in(rx)A$.
\medbreak
By the symmetry of \ref{tf16}, it suffices to do Cases 1--4. The proof of (a) with $\ell(x)=1$ is thus finished.

We will now prove (a) by induction on $\ell(x)$. For $\ell(x)\leq 1$, we have seen it before. Suppose $u\leq uv=x$, $u,v\in S_n-\{1\}$. Recall (\ref{tf34}): $\Inv(x)=\Inv(u)\amalg u\,\Inv(v)\,u^{-1}$ where $\amalg$ denotes disjoint union.

By the induction hypothesis applied to $(A,v)$ (instead of $(A,x)$), we have $\Inv(v)\subset(rv)A$. Hence
$(ru)\Inv(v)\subset (ru)(rv)A=(rx)A$. But the induction hypothesis for $(\Inv(v),u)$ implies that $(ru)\Inv(v)$ equals
the greatest half-permutation $B$ with $\Inv(u)\subset B\subset\Inv(u)\cup u\,\Inv(v)\,u^{-1}=\Inv(x)$. As $\Inv(x)$
is itself a half-permutation, we find $(ru)\Inv(v)=\Inv(x)$. We have thus shown:
\be \Inv(x)\subset(rx)A. \label{tf37} \ee
Applying the induction hypothesis to $(A,v)$, we see that  $(rv)A\subset\Inv(v)\cup vAv^{-1}$. Combining with the
induction hypothesis on $((rv)A,u)$, we find $(rx)A=(ru)(rv)A$ $\subset\Inv(u)\cup u(rv)Au^{-1}\subset\Inv(u)\cup
u(\Inv(v)\cup vAv^{-1})u^{-1}=(\Inv(u)\cup u\,\Inv(v)\,u^{-1})\cup uvAv^{-1}u^{-1}=\Inv(x)\cup xAx^{-1}$. We have
shown:
\be (rx)A\subset\Inv(x)\cup xAx^{-1}. \label{tf38} \ee
In view of (\ref{tf37}) and (\ref{tf38}), it remains to show that for any half-permutation $B$ with (\ref{tf21}) one has $B\subset (rx)A$. Let $B\in\HP$ have the property (\ref{tf21}). We have $\Inv(u)\subset\Inv(x)\subset B$, which shows
\be \Inv(u)\subset B. \label{tf40} \ee
Define $C=u^{-1}(B-\Inv(u))u$. By \ref{tf18} and (\ref{tf40}), we have
\be C\in \HP. \label{tf41} \ee
We have $\Inv(u)\amalg u\,\Inv(v)\,u^{-1}=\Inv(x)\subset B=\Inv(u)\amalg uCu^{-1}$, which shows
\be \Inv(v)\subset C. \label{tf42} \ee
We also have $\Inv(u)\amalg uCu^{-1}=B\subset\Inv(x)\cup xAx^{-1}=\Inv(u)\cup u\,\Inv(v)\,u^{-1}\cup uvAv^{-1}u^{-1}$, which shows
\be C\subset\Inv(v)\cup vAv^{-1}. \label{tf43}   \ee
Applying the induction hypothesis to $(A,v)$ and invoking (\ref{tf41}), (\ref{tf42}), (\ref{tf43}), one finds $C\subset(rv)A$. Hence $B=\Inv(u)\cup uCu^{-1}\subset\Inv(u)\cup u(rv)Au^{-1}$. Combining with (\ref{tf40}), we have $\Inv(u)\subset B\subset\Inv(u)\cup u(rv)Au^{-1}$. By the induction hypothesis applied to $((rv)A,u)$, it follows that $B\subset (ru)(rv)A=(rx)A$. This finishes the induction step and hence the proof of (a).

We turn to (b). By (a), $(rx)\Inv(y)$ 
is the greatest half-permutation $B$ with $\Inv(x)\subset B\subset \Inv(x)\cup x\,\Inv(y)\,x^{-1}=\Inv(xy)$, the last
identity being (\ref{tf34}). But $\Inv(xy)$ is itself a half-permutation. This proves that
$(rx)\Inv(y)=\Inv(xy)$.\enddemo

Recall the map $\,\GB=r\,L^{-1}\,\Pro\co\HP\ra\Omega$.

\proclaim{Lemma} \label{tf15} The map $\,\GB\co \HP\ra\Omega$ is $B_n^+$\/{\rm -}\/equivariant{\rm .}
 In formula{\rm ,} if $x\in B_n^+${\rm ,} $A\in\HP${\rm ,} $y=\GB(A)$ then $\GB(xA)=\LF(xy)${\rm .}\endproclaim

\demo{Proof} It suffices to give a proof for $\ell(x)=1$, so we will henceforth assume this is the case.
 Write $x=ru$, $y=rv$ ($u,v\in S_n$), and note $L(u)=\{u\}$. We know: \begin{itemize}
\item $\Pro((ru)A)$ equals the greatest $B\in\Inv(S_n)$ with $\{u\}\subset B\subset (ru)A$ (by Lemma \ref{tf14} and the observation $\Inv(S_n)\ni \{u\}\subset (ru)A$).
\item $(ru)A$ is the greatest $C\in\HP$ with $\{u\}\subset C\subset \{u\}\,\cup\, uAu$ (by\break
Lemma~\ref{tf20}(a)).\end{itemize} Combining these observations and recalling that $\Inv(S_n)\subset \HP$,
we immediately find that $\Pro((ru)A)$ is the greatest $B\in\Inv(S_n)$ with 
\be \{u\}\subset B\subset\{u\}\cup uAu. \label{tf23} \ee
Write $B=\Inv(uw)$, $w\in S_n$. Assume the left-hand inclusion of (\ref{tf23}) to hold:  $\{u\}\subset B$ or,
equivalently, $u\leq uw$. We have 
\begin{eqnarray}
\mbox{right-hand inclusion of }(\ref{tf23}) &\Longleftrightarrow& \Inv(uw)\subset\{u\}\cup uAu \nonumber \\
& \Longleftrightarrow & \{u\}\amalg u\,\Inv(w)\,u\subset \{u\}\cup uAu \nonumber \\
&\Longleftrightarrow& \Inv(w)\subset A \nonumber \\ &\stackrel{(\ref{tf14})}{\Longleftrightarrow}&
\Inv(w)\subset\Pro(A)=\Inv(v) \nonumber \\
&\Longleftrightarrow& w\leq v \nonumber \\
&\Longleftrightarrow& uw\leq r^{-1}\,\LF(xy), \nonumber \end{eqnarray}
the last equivalence following from the assumption that $u\leq uw$. The greatest $B$ satisfying these properties is given by $uw=r^{-1}\,\LF(xy)$. This shows $\Pro((ru)A)=\Inv\,r^{-1}\,\LF(xy)$ and the lemma follows.\enddemo

\section{Two more properties of  the representation} \label{secf}
\advance\eqcount by 16

Let $M_m(R)$ denote the algebra of size $m$ square matrices over $R$. We identify $M_m(R)$ with $\End(V)$.

 \proclaim{Theorem} \label{tf26} Suppose $R=\zz[q^{\pm1},t^{\pm1}]${\rm ,}
 the Laurent polynomial ring in two variables{\rm .} Let $x\in B_n${\rm ,}
 and consider the Laurent expansion of $\rho x$ with
respect to $t${\rm :}
$$ \rho x=\sum_{i=k}^\ell A_i(q)\,t^i,\ \ \ A_i\in M_m(\zz[q^{\pm1}]),\ \ \ A_k\neq0,\ \ \ A_\ell\neq0. $$
\begin{itemize}
\ritem{(a)} Then $\ell_\Omega(x)=\max(\ell-k,\ell,-k)${\rm .}
\ritem{(b)} If in addition $x\in B_n^+ -\Delta B_n^+${\rm ,}
 then $k=0$ and $\ell=\ell_\Omega(x)${\rm .}\end{itemize}
\endproclaim

{\it Proof}. First we prove (H1):  If $x\in B_n^+ -\Delta B_n^+$ then $k=0$. While we are primarily interested in
the case where $R$ is a Laurent polynomial ring, it obviously suffices to prove (H1) in the case where
$R=\rr[t^{\pm1}]$ and $q\in\rr$ with $0<q<1$. During the proof of (H1), we will assume this is the case. Since
$x\not\in\Delta B_n^+$ we have
$\LF(x)\neq\Delta$;   hence $C_{\LF(x)}\cap t\,V_1=\emptyset$. Choose any $v\in C_1$. Then $xv\in
C_{\LF(x)}$ whence $xv\not\in t\,V_1$.

Recall that $\rho B_n^+\subset M_m(\zz[q,q^{-1},t])$. 
Since $x\in B_n^+$, we thus have $k\geq 0$. Assume now $k>0$. Then all entries of $\rho x$ are in $t\,\rr[t]$;
hence $xv\in t\,V_1$. This is a contradiction. This proves that $k=0$; i.e., (H1) has been proved.

Next, we will show (H2):  If $x\in B_n^+ -\Delta B_n^+$ then $\ell=\ell_\Omega(x)$. Write $p=\ell_\Omega(x)$. Define $\Gamma\in\GL(V)$ by $\rho\Delta=t\Gamma$. By \ref{tf8} we have
\begin{eqnarray}
 T(q)\,\rho(x,q^{-1},t^{-1})\,T(q)^{-1}&=&\rho(\barr{x},q,t)=
\rho(\Delta^{-p},q,t)\,\rho(\Delta^p\barr{x},q,t)\label{tf27}\\
&=&\Gamma^{-p}t^{-p}\,\rho(\Delta^p\barr{x},q,t).\nonumber
\end{eqnarray}
 Note
that $T(q)$ does not involve $t$. Neither does $\Gamma$, by \ref{tf24}. We now compare the least exponents of $t$
occurring on either side of (\ref{tf27}). (For a matrix $A\in\GL(V)$, the least exponent of $t$ is by definition the
greatest integer $a$ such that $A\in t^a M_m(\zz[q,q^{-1},t])$.)  By Theorem~\ref{fa74}(d), we have
$\Delta^p\barr{x}\in B_n^+ -\Delta B_n^+$. Applying (H1) to $\Delta^p\barr{x}$ then shows that the least exponent of
$t$ on the right-hand side of (\ref{tf27}) equals $-p$. The least exponent on the left-hand side equals $-\ell$. It
follows that $\ell=p=\ell_\Omega(x)$. This finishes the proof of (H2), and hence of (b).

Finally, we prove (a). Recall the bijection $\zz\times(B_n^+ -\Delta B_n^+)\ra B_n$, $(a,y)\mapsto \Delta^ay$. Write $x=\Delta^ay$, $a\in\zz$, $y\in B_n^+ -\Delta B_n^+$, $\ell_\Omega(y)=b$. Then $k=a$, $\ell=a+b$ by (b) and \ref{tf24}. Using 
Theorem~\ref{fa74}(c), one finds $\ell_\Omega(x)=\ell_\Omega(\Delta^ay)=\max(a+b,b,-a)=\max(\ell,\ell-k,-k)$.
This proves (a).\hfill\qed\vglue6pt

An immediate consequence of Theorem \ref{tf26} is another proof of the faithfulness of $\rho\co B_n\ra\GL(V)$ (Theorem~\ref{tf7}). Indeed, if $x\in B_n$ is in the kernel of $\rho$, then in the notation of \ref{tf26}, we have $k=\ell=0$, whence $\ell_\Omega(x)=0$. It follows that $x=1$.

We return to our assumption $R=\rr[t^{\pm1}]$, $q\in\rr\subset R$,
$0<q<1$. Before proving our next theorem, we establish a simple lemma.
The results of Section~\ref{secb} (or see \cite{cha}) imply that any two positive braids $x,y$ have a greatest common lower bound, notation $x\wedge y$. For any two subsets $X,Y$ of some additive abelian group, we write $X+Y=\{x+y\mid x\in X,\ y\in Y\}$.

\proclaim{Lemma} \label{tf31} \begin{itemize}
\ritem{(a)} Let $s\in S$, $A\in\HP${\rm ,} $x=\GB(A)${\rm .} Then $rs\leq x\Leftrightarrow s\in A${\rm .}
\ritem{(b)} Let $x,y\in\Omega${\rm ,} $x\wedge y=1${\rm .} Then $C_x+C_y\subset C_1${\rm .}
 \end{itemize}\endproclaim

\demo{Proof} (a) We have $rs\leq x\Leftrightarrow s\leq r^{-1}x\Leftrightarrow \{s\}\subset \Inv(r^{-1}x)\Leftrightarrow \{s\}\subset\Pro(A)\Leftrightarrow \{s\}\subset A$. The last equivalence holds because $\Pro(A)$ is the greatest element of $\Inv(S_n)$ contained in $A$ by \ref{tf14}, and $\{s\}\in\Inv(S_n)$. This proves (a).

(b) Let $A\in \Brd^{-1}(x)$, $B\in \Brd^{-1}(y)$. 
We must show $D_A+D_B\subset C_1$. The intersection of any two half-permutations is again a half-permutation, so
$A\cap B\in\HP$. Note $D_A+D_B=D_{A\cap B}\subset C_z$ where $z=\Brd(A\cap B)$. We must therefore show $z=1$.
Suppose $z\neq 1$, say $s\in S$, $rs\leq z$. By (a), we have $s\in A\cap B$. By the other direction of (a) and the fact
that $s\in A$, we have $rs\leq x$. Similarly, $rs\leq y$; hence $rs\leq x\wedge y=1$. This contradiction
shows
$z=1$ and thus finishes the proof.\enddemo

\demo{Definition} We define a (total) ordering on $R=\rr[t^{\pm1}]$ as follows. Let $a\in R-\{0\}$, and write
$a=\sum_{i=k}^\ell a_it^i$, $a_i\in\rr$, $a_k\neq0$. Then the sign of $a$ is defined to be the sign of $a_k$. (This is
the only ordering of the ring $R$ which restricts to the usual ordering on $\rr$ and with $0<t<b$ for all positive real
numbers $b$.) We also define a map $\TP\co R\ra t^\zz\cup\{0\}$ (Trailing Power) which in the above notation takes
$a$ to $t^k$, and with $\TP(0)=0$.

\demo{Definition} We write $C$ instead of $C_1$. The union of all $xC$ (with $x\in B_n$) will be denoted by $U$.

Obviously, $C$ is closed under addition and scalar multiplication by elements of $\{a\in R\mid \TP(a)=1$ and $a>0\}$. We also have
\be aC=\TP(a)C\ \ \mbox{ for all }a\in R_{>0}. \label{tf44} \ee

The following theorem shows that $U$ has properties resembling those of convex cones in real vector spaces, and moreover relates the greedy form with line segments in $U$ defined over $\rr[t^{\pm1}]$.

 \proclaim{Theorem} \label{tf28} \begin{itemize}
\ritem{(a)} $\Delta C=tC${\rm .}
\ritem{(b)} The $xC$ {\rm (}\/with $x\in B_n$\/{\rm )} are disjoint{\rm .}
\ritem{(b)} Let $(y_1,\ldots,y_k)\in\Omega^k$ be greedy\/{\rm ;} i.e.{\rm ,} $\LF(y_iy_{i+1})=y_i$ $(1\leq i<k)${\rm .}
 Let
$x_0,\ldots,x_k\in B_n$ be such that $x_i=x_{i-1}y_i$ $(1\leq i\leq k).$ Then 
$$ t^ix_0C+x_kC\subset \left\{ \begin{array}{llr}
t^i {x_0}C,\ \ & i\leq 0 & {  \ \ \ \ ({\rm c1})}; \\
{x_i}C, & 0\leq i\leq k & {  \ ({\rm c2})}; \\
{x_k}C, & k\leq i & {  \ ({\rm c3})}. \\ \end{array} \right. $$
\ritem{(d)} Let $(\tilde{y}_1,\ldots,\tilde{y}_k)$ be a Thurston normal form\/{\rm ;} i.e.\/{\rm ,}\/
 there are greedy\break
$(u_1,\ldots,u_s)$\/{\rm ,}\/
$(v_1,\ldots,v_t)$ with $(u_s^{-1},\ldots,u_1^{-1},v_1,\ldots,v_t)=(\tilde{y}_1,\ldots,\tilde{y}_k)$\/{\rm ,}\/
 and $u_s,v_t\neq
1$\/{\rm ,}\/ and there is no $w\in B_n^+ -\{1\}$ such that $\{u_1,v_1\}\subset wB_n^+$\/{\rm .}\/
 Let $\tilde{x}_0,\ldots,\tilde{x}_k\in
B_n$ be such that $\tilde{x}_i=\tilde{x}_{i-1}\tilde{y}_i$ $(1\leq i\leq k).$ Then 
$$ \frac{t^i\tilde{x}_0C+\tilde{x}_kC}{t^i+1}\subset \left\{ \begin{array}{llr}
{\tilde{x}_0}C,\ \ & i\leq -s; \\
{\tilde{x}_{i+s}}C, & -s\leq i\leq t; \\
{\tilde{x}_k}C, & t\leq i. \\ \end{array} \right. $$
\ritem{(e)} The set $U$ is closed under addition and scalar multiplication by positive elements of
$R${\rm .}\end{itemize}

\endproclaim
{\it Proof}. (a) This follows from \ref{tf24}: 
 $\Delta\,x_{n+1-j,n+1-i}=tq^{i+j-1}\,x_{ij}$ whenever $1\leq i<j\leq n$, and the involution of \ref{tf16}.
\smallbreak
(b) Let $x\in B_n$, $x\neq 1$. We must show that $C$ and
 $xC$ are disjoint. Write $x=y\Delta^k$, $y\in B_n^+ -\Delta B_n^+$, $k\in \zz$. Note:
\be C\subset V_1-tV_1. \label{tf29} \ee
Similarly, we have $yC=yC_1\subset C_{\LF(y)}\subset V_1-t V_1$ (because $y\not\in\Delta B_n^+$) whence
\be xC=y\Delta^k C=t^kyC\subset t^kV_1-t^{k+1}V_1. \label{tf30} \ee
Suppose now $k>0$. Then (\ref{tf30}) shows $xC\subset t^k V_1\subset tV_1$. Combining with (\ref{tf29}), one finds that $C$ and $xC$ are disjoint. Next, suppose $k<0$. Then (\ref{tf29}) shows $C\subset V_1\subset t^{k+1} V_1$, which cannot meet $xC$ by (\ref{tf30}). It remains to consider the case $k=0$, i.e., $y=x$. Suppose $C\cap xC\neq\emptyset$. We have $xC\subset C_{\LF(x)}$ and $C=C_1$, so that $C_1\cap C_{\LF(x)}\neq\emptyset$. Since all $C_z$ are disjoint by \ref{tf12}(a), it follows that $\LF(x)=1$, whence $x=1$. This finishes the proof in the case $k=0$, and thereby proves (b).
\smallbreak
(c1) Let $i\leq0$. Then $t^ix_0C+x_kC=t^ix_0(C+t^{-i}(y_1\cdots y_k)C)\subset t^ix_0(C+t^{-i}V_2)\subset t^ix_0(C+V_2)=t^ix_0C$. This proves (c1).
\smallbreak
(c2) First, we consider the case $i=1$. Note that $y_1^{-1}\Delta\in\Omega$. Moreover, the fact that $(y_1,y_2)$ is greedy (i.e., $\LF(y_1y_2)=y_1$) is equivalent to $(y_1^{-1}\Delta)\wedge y_2=1$. Using (a), we find
$$ x_1^{-1}(t{x_0}C+{x_k}C)\stackrel{{\rm (a)}}{=}
(y_1^{-1}\Delta)C_1+(y_2\cdots y_k)C_1\subset C_{y_1^{-1}\Delta}+C_{y_2}\subset C_1=C.  $$
Here, the first inclusion follows from \ref{tf12}(b), and the second inclusion from \ref{tf31}(b). This finishes the proof in the case $i=1$. We now give a proof of (c2) by induction on $i$. For $i=0$, it follows from (c1). The induction step is shown as follows. Notice $C=C+tC$. Hence
\begin{eqnarray*}
 t^ix_0C+x_kC&=&t^ix_0C+(tx_kC+x_kC)\\
&=&t(t^{i-1}x_0C+x_kC)+x_kC\subset tx_{i-1}C+x_kC\subset x_iC.   
\end{eqnarray*}
Here, the first inclusion is the induction hypothesis and the second inclusion is a shifted version of the $i=1$ case. This finishes the proof of (c2).
\smallbreak
(c3) Notice that in (c2), it is not excluded that some $y_i$ is 1. By extending the sequence $(y_1,\ldots,y_k)$ in (c1) far enough to the right by ones, one can assume some new $k$ to be at least $i$. Then (c2) applies and shows (c3).
\smallbreak
(d) Define $x_i\in B_n$ by
$$ x_i= \left\{ \begin{array}{ll}
\tilde{x}_i\Delta^{i-s}, & 0\leq i\leq s; \\
\tilde{x}_i, & s\leq i\leq k. \\ \end{array} \right. $$
Define $y_i\in \Omega$ by $x_i=x_{i-1}y_i$ ($1\leq i\leq k$). Then $x_i$ and $y_i$ are as in (c). Observe:
$$ x_iC=t^{i-s}\tilde{x}_iC \ \ \  (0\leq i\leq s). $$
We have
$$ \frac{t^i\tilde{x}_0C+\tilde{x}_kC}{t^i+1}=\frac{t^{i+s}x_0C+x_kC}{t^i+1}. $$
The inclusions in the sequel will be consequences of (c). If $i\leq -s$ then
$$ \frac{t^{i+s}x_0C+x_kC}{t^i+1} \subset \frac{t^{i+s}x_0C}{t^i+1}=t^sx_0C=\tilde{x}_0C. $$
If $-s\leq i\leq 0$ then
$$ \frac{t^{i+s}x_0C+x_kC}{t^i+1} \subset \frac{x_{i+s}C}{t^i+1}= \frac{t^i\tilde{x}_{i+s}C}{t^i+1}=\tilde{x}_{i+s}C. $$
If $0\leq i\leq t$ then
$$ \frac{t^{i+s}x_0C+x_kC}{t^i+1} \subset \frac{x_{i+s}C}{t^i+1}= x_{i+s}C=\tilde{x}_{i+s}C. $$
Finally, if $t\leq i$ then
$$ \frac{t^{i+s}x_0C+x_kC}{t^i+1} \subset \frac{x_kC}{t^i+1}=x_kC=\tilde{x}_kC. $$
This proves (d).
\smallbreak
Part (e) is an easy consequence of (a), (c) and (\ref{tf44}).\hfill\qed

\end{document}